\documentclass[10pt]{article}

\usepackage{amsmath}
\usepackage{amsfonts}
\usepackage{amssymb}

\newtheorem{proposition}{Proposition}[section]
\newtheorem{theorem}{Theorem}[section]
\newtheorem{lemma}[theorem]{Lemma}
\newtheorem{corollary}[theorem]{Corollary}
\newtheorem{remark}[theorem]{Remark}

\def\phi{{\varphi}}

\DeclareSymbolFont{AMSb}{U}{msb}{m}{n}
\DeclareMathSymbol{\N}{\mathbin}{AMSb}{"4E}
\DeclareMathSymbol{\Z}{\mathbin}{AMSb}{"5A}
\DeclareMathSymbol{\R}{\mathbin}{AMSb}{"52}
\DeclareMathSymbol{\Q}{\mathbin}{AMSb}{"51}
\DeclareMathSymbol{\I}{\mathbin}{AMSb}{"49}
\DeclareMathSymbol{\C}{\mathbin}{AMSb}{"43}

\begin{document}

\title{On the best possible remaining term in the Hardy inequality}
\author{ Nassif  Ghoussoub\thanks{Partially supported by a grant
from the Natural Sciences and Engineering Research Council of Canada.  } \quad  and \quad Amir Moradifam 
\\
\small Department of Mathematics,
\small University of British Columbia, \\
\small Vancouver BC Canada V6T 1Z2 \\
\small {\tt nassif@math.ubc.ca} \\
\small {\tt a.moradi@math.ubc.ca}
\\
}
\maketitle

\section{Abstract} 
We give a necessary and sufficient condition on  a radially symmetric potential $V$ on  $\Omega$ that makes it  an admissible candidate for an  improved Hardy inequality of the following form:
\begin{equation}\label{gen-hardy.0}
\hbox{$\int_{\Omega}|\nabla u |^{2}dx - ( \frac{n-2}{2})^{2} \int_{\Omega}\frac{|u|^{2}}{|x|^{2}}dx\geq c\int_{\Omega} V(|x|)|u|^{2}dx$ \quad for all $u \in H^{1}_{0}(\Omega)$.}
\end{equation}
A characterization of the best possible constant $c(V)$ is also given. This result yields easily the improved Hardy's inequalities of Brezis-V\'{a}zquez \cite{BV},  Adimurthi et al. \cite{ACR}, and Filippas-Tertikas \cite{FT}, as well as the corresponding best constants. Our approach clarifies the issue behind the lack of an optimal improvement, while yielding other interesting ``dual"  inequalities. Another consequence is the following substantial sharpening of known integrability criteria: If a  positive radial function $V$  satisfies 
$\hbox{$\liminf_{r\rightarrow 0} \ln(r)\int^{r}_{0} sV(s)ds>-\infty$, } $
then there exists $\rho:=\rho(\Omega)>0$ such that the improved Hardy inequality (\ref{gen-hardy.0}) holds for the scaled potential $V_\rho(x)=V(\frac{|x|}{\rho})$.
 On the other hand, if 
$\hbox{$\lim_{r\rightarrow 0} \ln(r)\int^{r}_{0} sV(s)ds=-\infty$, } $
then there is no $\rho>0$ for which (\ref{gen-hardy.0}) holds for $V_\rho$.  This shows for example,  that $V(x)=\frac{1}{x^\alpha}$ is an admissible potential for an improved Hardy inequality when $\alpha <2$, while it is not so for $\alpha \geq 2$. All these results have immediate applications to the corresponding  Schr\"odinger equations. Analogous criteria for inequalities involving the  bi-Laplacian will be developed in a forthcoming paper \cite{GM2}.

\section{Introduction}
Let $\Omega$ be a bounded domain in $R^{n}$, $n \geq 3$, with $0 \in \Omega$. The classical Hardy inequality asserts that 
\begin{equation}\label{cl-hardy}
\hbox{$\int_{\Omega}|\nabla u |^{2}dx \geq ( \frac{n-2}{2})^{2} \int_{\Omega}\frac{|u|^2}{|x|^{2}}dx$ \quad  for all $u \in H^{1}_{0}(\Omega)$.}
\end{equation}
This inequality and its various improvements  are used in many contexts, such as in the study of the stability of solutions of semi-linear elliptic and parabolic equations \cite{BV, CM1, V}, in the analysis of the asymptotic behavior of the heat equation with singular potentials \cite{CM2, VZ}, as well as in the study of the stability of eigenvalues in elliptic problems such as Schr\"odinger operators \cite{D,FHT}.

Now it is well known that $( \frac{n-2}{2})^{2} $ is the best constant for inequality (\ref{cl-hardy}), and that this constant is however not attained in $H^{1}_{0}(\Omega)$. So, one could anticipate  
improving this inequality by adding a non-negative correction term to the right hand side of (\ref{cl-hardy}) and   
indeed, several sharpened Hardy inequalities have been established in recent years \cite{BM, BMS, FT, FHT, VZ}, mostly  triggered  
by the following improvement of Brezis and V\'{a}zquez \cite{BV}. 
\begin{equation}\label{bv-hardy}
\hbox{$\int_{\Omega}|\nabla u |^{2}dx \geq ( \frac{n-2}{2})^{2} \int_{\Omega}\frac{|u|^2}{|x|^{2}}dx+\lambda_{\Omega}\int_{\Omega}|u|^{2}dx$ \quad for every $u \in H^{1}_{0}(\Omega)$.}
\end{equation}
 The constant $\lambda_{\Omega}$ in (\ref{bv-hardy}) is given by
\begin{equation}\label{constant}
\lambda_{\Omega}=z^{2}_{0}\omega^{2/n}_{n}|\Omega|^{-\frac{2}{n}},
\end{equation}
where $\omega_{n}$ and $|\Omega|$ denote the volume of the unit ball and $\Omega$ respectively, and $z_{0}$ is the first zero of the bessel function $J_{0}(z)$. Moreover,  $\lambda_{\Omega}$ is optimal when $\Omega$ is a ball, but  is --again-- not achieved in $H^{1}_{0}(\Omega)$. This led to one of the open problems mentioned in \cite{BV} (Problem 2), which is whether the two terms on the RHS of inequality (\ref{bv-hardy}) (i.e., the coefficients of $|u|^2$) are just the first two terms of an infinite series of correcting terms. 

This question was addressed by several authors. In particular, Adimurthi et all \cite{ACR} proved that for every integer $k$, there exists a constant $c$ depending on $n$, $k$ and $\Omega$ such that
\begin{equation}\label{ar-hardy}
\hbox{$\int_{\Omega}|\nabla u |^{2}dx \geq ( \frac{n-2}{2})^{2} \int_{\Omega}\frac{|u |^{2}}{|x|^{2}}dx+c\sum^{k}_{j=1}\int_{\Omega}\frac{|u |^{2}}{|x|^{2}}\big( \prod^{j}_{i=1}log^{(i)}\frac{\rho}{|x|}\big)^{-2}dx$ \quad for $u \in H^{1}_{0}(\Omega)$,}
\end{equation}
where $\rho=(\sup_{x \in \Omega}|x|)( e^{e^{e^{.^{.^{e(k-times)}}}}} )$. Here we have used the notations 
$log^{(1)}(.)=log(.)$ and  $log^{(k)}(.)=log(log^{(k-1)}(.))$ for $k\geq 2$. 

Also motivated by the question of Brezis and V\'{a}zquez, Filippas and Tertikas  proved in \cite{FT} that the inequality can be repeatedly improved by adding to the right hand side specific potentials which lead to an infinite series expansion of Hardy's inequality. More precisely, by defining iteratively the following functions, 
\begin{equation*}\label{x-def}
X_{1}(t)=(1-\log(t))^{-1}, \quad  X_{k}(t)=X_{1}(X_{k-1}(t)) \ \ \ \ k=2,3, ... ,
\end{equation*}
they prove that for any $D\geq \sup_{x \in \Omega}|x|$, the following inequality holds for any $u \in H_{0}^{1}(\Omega)$:
\begin{equation}\label{ft-hardy}
\int_{\Omega}|\nabla u |^{2}dx \geq ( \frac{n-2}{2})^{2} \int_{\Omega}\frac{|u|^2}{|x|^{2}}dx+\frac{1}{4}\sum^{\infty}_{i=1}\int_{\Omega}\frac{1}{|x|^{2}}X^{2}_{1}(\frac{|x|}{D})X^{2}_{2}(\frac{|x|}{D})...X^{2}_{i}(\frac{|x|}{D})|u|^{2}dx.
\end{equation}
Moreover, they proved that the constant $\frac{1}{4}$ is the best constant for the corresponding $k-$improved Hardy inequality which is again not attained in $H_{0}^{1}(\Omega)$.

In this paper,  we  show that all the above results --and more--  follow from a  specific characterization of those potentials $V$ that yield an improved Hardy inequality.
Here is our main result.

\begin{theorem} \label{main} Let $V$ be a radial function on a smooth bounded radial domain $\Omega$ in $\R^n$ with radius $R$, in such a way that $V(x)= v(|x|)$ for some non-negative  function $v$ on $(0,R)$. The following properties are then equivalent:
\begin{enumerate}
\item The ordinary differential equation 
\[
\hbox{ $({\rm D}_V)$  \quad \quad \quad \quad \quad  \quad \quad \quad \quad \quad \quad \quad \quad \quad \quad $y''(r)+\frac{y'(r)}{r}+v(r)y(r)=0$  \quad \quad \quad \quad \quad  \quad \quad \quad \quad \quad  \quad \quad \quad \quad \quad \quad \quad \quad \quad \quad}
\]
 has a positive solution on the interval $(0, R)$.

\item The  following improved Hardy inequality holds
\begin{equation*}\label{gen-hardy}
\hbox{   $({\rm H}_{V})$ \quad \quad \quad \quad \quad \quad \quad \quad $\int_{\Omega}|\nabla u |^{2}dx - ( \frac{n-2}{2})^{2} \int_{\Omega}\frac{|u|^{2}}{|x|^{2}}dx\geq  \int_{\Omega} V(|x|)|u|^{2}dx$ \quad for $u \in H^{1}_{0}(\Omega)$. \quad \quad \quad \quad \quad \quad \quad  }
\end{equation*}
\end{enumerate} 
Moreover, the best constant \hbox{$c(V):=\sup\big\{c; \, \,  ({\rm H}_{cV})$ holds$\big\}$} can be characterized by the formula
\begin{equation} \label{best.constant}
\hbox{$c(V)=\sup\big\{c; \,  y''(r)+\frac{y'(r)}{r}+cv(r)y(r)=0$ has a positive solution on the interval $\big(0, R \big)\big\}$.}
\end{equation}
\end{theorem} 
We note that the implication 1) implies 2) holds for any smooth bounded domain $\Omega$ in $\R^n$ containing $0$, provided $v(r)+(\frac{n-2}{2})^{2}\frac{1}{r^{2}}$  is non-increasing on $(0,\sup_{x \in \Omega}|x|)$ and $R$ is the radius of the ball which has the same volume as $\Omega$ (i.e. $R=(\frac{|\Omega|}{\omega_{n}})^{\frac{1}{n}})$.

It is  therefore clear from the above discussion that in order to find what potentials are candidates for an improved Hardy inequality, one needs to investigate the ordinary differential equation
$y''+\frac{y'}{r}+v(r)y(r)=0.$
 We shall see that the results of Brezis-V\'{a}zquez,  Adimurthi et al, and Filippas-Tertikas mentioned above can be easily deduced by simply checking that the potentials $V$ they consider, correspond to equations  $({\rm D}_V)$ where an explicit positive solution can be found.

Our approach turned out to be also useful for determining the best constants in the above mentioned improvements. Indeed, the case when $V\equiv 1$ will follow  immediately from Theorem \ref{main}.  A slightly more involved reasoning -- but also based of the above characterization --  will allow us to find the best constant in the improvement of Adimurthi et al, and to recover the best one established by Filippas-Tertikas. 

Since the existence of positive solutions for ODEs of the form  $({\rm D}_V)$ is closely related to the oscillatory properties of second order equations of the form $z''(s)+ a(s)z(s)=0$, Theorem \ref{main} also  allows for the use of the extensive literature on the oscillatory properties of  such equations to deduce various interesting results such as the following corollary.  

\begin{corollary} \label{osc} Let $V$ be a positive radial function on a smooth bounded radial domain $\Omega$ in $\R^n$. 
\begin{enumerate}
\item If 
$\hbox{$\liminf_{r\rightarrow 0} \ln(r)\int^{r}_{0} sV(s)ds>-\infty$, } $
then there exists $\alpha:=\alpha(\Omega)>0$ such that an improved Hardy inequality $({\rm H}_{V_\alpha})$ holds for the scaled potential $V_\alpha(x):=\alpha^2V(\alpha x)$.  
\item  If
$\hbox{$\lim_{r\rightarrow 0} \ln(r)\int^{r}_{0} sV(s)ds=-\infty$, } $
then there are no $\beta, c>0$,  for which $({\rm H}_{V_{\beta,c}})$ holds with $V_{\beta,c}=cV(\beta x)$.
\end{enumerate}
\end{corollary}
The following is a consequence of the two results above. 

\begin{corollary} For any $\alpha<2$,  inequality $({\rm H}_{cV})$ holds on a bounded domain $\Omega$ 
for $V_{\alpha}(x)=\frac{1}{|x|^{\alpha}}$ and  some $c>0$. Moreover, the
best constant $c(\frac{1}{|x|^{\alpha}})$ is equal to the largest 
$c$ such that the equation
\begin{equation*}
y''(r)+\frac{1}{r}y'(r)+c\frac{1}{|x|^{\alpha}}=0,
\end{equation*}
has a positive solution on $(0,R)$, where $R$ is the radius of the ball
wich has the same volume as $\Omega$. Moreover, if $\alpha \geq 2$ inequality $({\rm H}_{V})$ does not hold for $V_{\alpha,c}(x)=c\frac{1}{|x|^{\alpha}}$ for any $c>0$. 
\end{corollary}
Note that the above corollary gives another proof of the fact that $(\frac{n-2}{2})^{2}$ is the best constant for the classical Hardy inequality.

Define now the class 
\[\hbox{$A_{\Omega}=\{v:R \rightarrow \R^+ ; $v is non-increasing on $(0,\sup_{x \in \Omega}|x|)$, $D_{v}$ has a positive solution on $(0,(\frac{|\Omega|}{\omega_{n}})^{\frac{1}{n}})$\}. }
\]
An immediate application of Theorem \ref{main} coupled with H\"older's inequality gives the following duality statement, which should  be compared to  inequalities  dual to those of Sobolev, recently obtained via the theory of mass transport  \cite{AGK, CNV}.

\begin{corollary} \label{dual} Suppose that $\Omega$ is a smooth bounded domain in $R^{n}$ containing $0$. Then for any $0<p\leq 2$, we have 
\begin{equation}
\inf \left\{\int_{\Omega}|\nabla u |^{2}dx - ( \frac{n-2}{2})^{2} \int_{\Omega}\frac{|u|^{2}}{|x|^{2}}dx; \, u \in H_{0}^{1}(\Omega), ||u||_p=1\right\}
\geq \sup \left\{\frac{1}{||V^{-1}(|x|)||_{L^{\frac{p}{p-2}}(\Omega)}.} ; \, V\in  A_{\Omega}\right\}.
\end{equation}
\end{corollary}

Finally, consider the following classes of radial potentials:
\begin{equation}
X=\{V:\Omega \rightarrow \R^+; \, V\in L^\infty_{\rm loc}(\Omega\setminus \{0\}), \ \ \liminf_{r\rightarrow 0} \ln(r)\int^{r}_{0} sV(s)ds>-\infty\},
\end{equation}
and
\begin{equation}
Y=\{V:\Omega \rightarrow \R^+;\, V\in L^\infty_{\rm loc}(\Omega\setminus \{0\}), \ \ \lim_{r\rightarrow 0} \ln(r)\int^{r}_{0} sV(s)ds=-\infty\}.
\end{equation}
 
For any $0<\mu<\mu_n:=\frac{(n-2)}{2})^{2}$ we consider the following weighted eigenvalue problem,
\begin{equation}
(E_{V, \mu}) \quad  \quad  \quad  \quad  \quad \quad  \quad  \quad  \quad  \quad \left\{ \begin{array}{lcl}
-\Delta u-\frac{\mu}{|x|^{2}}u&=&\lambda Vu \ \ {\rm in} \ \  \Omega, \\
\hfill u&=&0 \ \  \ \quad {\rm on} \ \ \Omega.
\end{array}
\right.
\quad  \quad  \quad  \quad  \quad \quad  \quad  \quad  \quad  \quad
\end{equation}
Our results above combine with standard arguments to yield the following. 
\begin{corollary} \label{eigenvalue} For any $0<\mu<\mu_n$, and $V:\Omega \rightarrow \R^+$ with $ V\in L^\infty_{\rm loc}(\Omega\setminus \{0\})$ and $ \lim_{|x|\rightarrow 0}|x|^{2}V(x)=0$, the  problem $(E_{V, \mu})$ admits a positive weak solution $u_\mu\in H^1_0(\Omega)$ corresponding to the first eigenvalue $\lambda=\lambda_\mu^1(V)$. Moreover, by letting $\lambda_1 (V)=\lim_{\mu \uparrow \mu_n}\lambda_\mu^1(V)$, we have
\begin{itemize}
\item If $V\in X$, then there exists $c>o$ such that $\lambda_1 (V_c)>0$.
\item If $V\in Y$, then $\lambda_1 (V_c)=0$ for all $c>0$,
\end{itemize}
where $V_{c}(x):=V(cx)$.
\end{corollary} 
 
 \section{Two dimensional inequalities}
In this section,  we start by establishing the following improvements of ``two-dimensional"  Poincar\'{e} and Poincar\'{e}-Wirtinger inequalities. 

\begin{theorem}\label{2dim-theorem}
Let $a<b$, $k$ is a differentiable function on $(a,b)$, and $\varphi$ be a strictly positive real valued differentiable function on $(a,b)$. Then,  every $h \in C^{1}([a,b])$ with
\begin{equation}\label{2dim-main-con}
-\infty<\lim_{r\rightarrow a}k(r)|h(r)|^{2}\frac{\varphi'(r)}{\varphi(r)}=\lim_{r\rightarrow b}k(r)|h(r)|^{2}\frac{\varphi'(r)}{\varphi(r)}<\infty,
\end{equation}
 satisfies the following inequality:
\begin{equation}\label{2dim-main-in}
\int^{b}_{a}|h'(r)|^{2}k(r) dr \geq \int^{b}_{a} -|h(r)|^{2}(\frac{k'(r)\varphi'(r)+k(r)\varphi''(r)}{\varphi(r)})dr.
\end{equation}
Moreover, assuming (\ref{2dim-main-con}),  the 
equality holds if and only if $h(r)=\varphi(r)$ for all $r \in (a,b)$.
\end{theorem}
Proof. Define $\psi(r)=h(r)/\varphi(r)$, $r \in [a,b]$. Then
\begin{eqnarray*}
\int^{b}_{a}|h'(r)|^{2}k(r) dr&=&\int^{b}_{a} |\psi (r)|^{2}|\varphi'(r)|^{2}k(r)dr+\int^{b}_{a} 2\varphi(r)\varphi'(r)\psi(r)\psi'(r)k(r)dr +\int^{b}_{a} |\varphi (r)|^2|\psi' (r)|^{2}k(r)dr \\
&=&\int^{b}_{a} |\psi (r)|^{2}|\varphi'(r)|^{2}k(r)dr-\int^{b}_{a}|\psi (r)|^{2}( k\varphi \varphi')'(r)dr+\int^{b}_{a} |\varphi (r)|^2|\psi' (r)|^{2}k(r)dr \\
&=&\int^{b}_{a} |\psi (r)|^{2}(|\varphi'(r)|^{2}k(r)-( k\varphi \varphi')'(r)dr+\int^{b}_{a} |\varphi (r)|^2|\psi' (r)|^{2}k(r)dr. 
\end{eqnarray*}
Hence, we have
\begin{eqnarray}
\int^{b}_{a}|h'(r)|^{2}k(r) dr&=&\int^{b}_{a} -|h(r)|^2(\frac{k'(r)\varphi'(r)+k(r)\varphi''(r)}{\varphi})dr+\int^{b}_{a} |\varphi (r)|^2|\psi' (r)|^{2}k(r)dr \\
&\geq& \int^{b}_{a} -|h(r)|^2(\frac{k'(r)\varphi'(r)+k(r)\varphi''(r)}{\varphi(r)})dr. 
\end{eqnarray}
Hence (\ref{2dim-main-in}) holds. Note that the last inequality is obviously an idendity if and only if $h(r)=\varphi(r)$ for all $r \in(a,b)$. The proof is complete. $\square$ \\

By applying Theorem \ref{2dim-theorem} to the weight  $k(r)=r$, we obtain the following generalization of the $2$-dimensional  Poincar\'{e} inequality.
\begin{corollary}\label{poinc}{\rm (Generalized 2-dimensional Poincar\'{e} inequality)} Let $0\leq a<b$ and $\varphi$ be a strictly positive real valued differentiable function on $(a,b)$. Then every $h \in C^{1}([a,b])$ with
\begin{equation}\label{poinc-con}
-\infty<\lim_{r\rightarrow a}r|h(r)|^{2}\frac{\varphi'(r)}{\varphi(r)}=\lim_{r\rightarrow b}r|h(r)|^{2}\frac{\varphi'(r)}{\varphi(r)}<\infty,
\end{equation}
satisfies the following inequality:
\begin{equation}
\int^{b}_{a}|h'(r)|^{2}r dr \geq \int^{b}_{a} -|h(r)|^{2}(\frac{\varphi'(r)+r\varphi''(r)}{\varphi(r)})dr.
\end{equation}
Moreover, under the assumption (\ref{poinc-con}) the equality holds if and only if $h(r)=\varphi(r)$ for all $r \in (a,b)$.
\end{corollary}
By applying Theorem \ref{2dim-theorem} to the weight  $k(r)=1$, we obtain the following generalization of the $2$-dimensional  Poincar\'{e}-Wirtinger inequality.

\begin{corollary}\label{poinc-wir}{\rm (Generalized Poincar\'{e}-Wirtinger inequality)} Let $a<b$ and $\varphi$ be a strictly positive real valued differentiable function on $(a,b)$. Then, every $h \in C^{1}([a,b])$ with
\begin{equation}\label{poinc-wir-con}
-\infty<\lim_{r\rightarrow a}|h(r)|^{2}\frac{\varphi'(r)}{\varphi(r)}=\lim_{r\rightarrow b}|h(r)|^{2}\frac{\varphi'(r)}{\varphi(r)}<\infty,
\end{equation}
satisfies the following inequality:
\begin{equation}
\int^{b}_{a}|h'(r)|^{2} dr \geq \int^{b}_{a} -|h(r)|^{2}\frac{\varphi''(r)}{\varphi(r)}dr.
\end{equation}
Moreover, under assumption (\ref{poinc-wir-con}),  the equality holds if and only if $h(r)=\varphi(r)$ for all $r \in (a,b)$.
\end{corollary}

\begin{remark}
Note that all of inequalities presented in the above theorems hold when we replace the condtions (\ref{2dim-main-con}), (\ref{poinc-con}),
and (\ref{poinc-wir-con}) with the following weaker conditions
\begin{eqnarray*}
\liminf_{r\rightarrow b}k(r)|h(r)|^{2}\frac{\varphi'(r)}{\varphi(r)}&\geq& \limsup_{r\rightarrow a}k(r)|h(r)|^{2}\frac{\varphi'(r)}{\varphi(r)},\\
\liminf_{r\rightarrow b}r|h(r)|^{2}\frac{\varphi'(r)}{\varphi(r)}&\geq& \limsup_{r\rightarrow a}r|h(r)|^{2}\frac{\varphi'(r)}{\varphi(r)},\\
\liminf_{r\rightarrow b}|h(r)|^{2}\frac{\varphi'(r)}{\varphi(r)}&\geq& \limsup_{r\rightarrow a}|h(r)|^{2}\frac{\varphi'(r)}{\varphi(r)},
\end{eqnarray*}
respectively, provided both sides in the above inequalities are not equal to $-\infty$ or $+\infty$.
\end{remark}

 \section{Proof of  the main theorem} 

 We start with the sufficient condition of Theorem 2.1 by establishing the following.
 
  \begin{proposition}\label{main-pro}  (Improved Hardy Inequality) Let  $\Omega$ be a bounded  smooth domain in $R^{n}$  with $0 \in \Omega$, and set $R=(|\Omega|/\omega_{n})^{1/n}$. Suppose $V$ is a radially symmetric function on $\Omega$ and  $\varphi$ is a $C^2$-function on $(0, R)$ such that
 \begin{equation}
\hbox{$0\leq V(|x|)\leq-\frac{\varphi'(|x|)+r\varphi''(|x|)}{|x|\varphi(|x|)}$ for all $x\in \Omega, \, 0<|x|<R$, }
\end{equation}
 \begin{equation}\label{boundary}
\hbox{$\liminf\limits_{r\rightarrow 0}r\frac{\varphi'(r)}{\varphi(r)}\geq 0$ \quad  and  \quad $\limsup\limits_{r\rightarrow R}\frac{\varphi'(r)}{\varphi(r)}<\infty$,}
\end{equation} 
 \begin{equation}
\hbox{ $(\frac{n-2}{2})^{2}\frac{1}{|x|^{2}}+V(|x|)$  is a decreasing function of $|x|$.} 
 \end{equation}
 Then for any $u \in H^{1}_{0}(\Omega)$, we have 
\begin{equation}\label{main-in}
\int_{\Omega}|\nabla u |^{2}dx \geq ( \frac{n-2}{2})^{2} \int_{\Omega}\frac{|u|^2}{|x|^{2}}dx+\int_{\Omega}V(|x|)|u|^{2}dx.
\end{equation}
Moreover, if $\lim_{r\rightarrow 0}r\varphi(r)\varphi'(r)=\lim_{r\rightarrow R}\varphi(r)\varphi'(r)=0$, then equality holds if and only if $u$ is a radial function on $\Omega$ such that $u(x)=\varphi(|x|)$ for all $x \in \Omega$.

\end{proposition}


\noindent {\bf Proof:} We  first prove the inequality for smooth radial positive functions on the ball $\Omega=B_R$. For such $u \in C^{2}_{0}(B_R)$,  we define
\[v(r)=u(r)r^{(n-2)/2}, \ \ r=|x|.\]
 In view of Corollary \ref{poinc}, we can write
\begin{eqnarray*}
\int_{\Omega}|\nabla u(x) |^{2}dx -( \frac{n-2}{2})^{2} \int_{\Omega}\frac{u^{2}(x)}{|x|^{2}}dx&=&\omega_{n}\int^{R}_{0}|\frac{n-2}{2}r^{-n/2}v(r)-r^{1-n/2}v'(r)|^{2}r^{n-1}dr\\
&-&( \frac{n-2}{2})^{2}\omega_{n}\int^{R}_{0}\frac{v^{2}(r)}{r}dr \\
&=&\omega_{n}( \frac{n-2}{2})^{2}\int^{R}_{0}v^{2}((1-\frac{2v'(r)r}{(n-2)v(r)})^{2}-1)\frac{dr}{r}\\
&=&\omega_{n}\int^{R}_{0} (v'(r))^{2}r-\omega_{n}( \frac{n-2}{2})\int^{R}_{0}v(r)v'(r)dr \\
&=&\omega_{n}\int^{R}_{0} (v'(r))^{2}r \\
&\geq & \omega_{n} \int^{R}_{0} -v^{2}(r)(\frac{\varphi'(r)+r\varphi''(r)}{\varphi(r)})dr \\
&=&\omega_{n} \int^{R}_{0} -u^{2}(r)(\frac{\varphi'(r)+r\varphi''(r)}{\varphi(r)})r^{n-2}dr \\
&=&-\int_{\Omega}u^{2}(x)(\frac{\varphi'(|x|)+|x|\varphi''(|x|)}{|x|\varphi(|x|)})dx.
\end{eqnarray*}
Hence, the inequality (\ref{main-in}) holds for radial smooth positive functions. By density arguments, inequality (\ref{main-in}) is valid for any $u \in H^{1}_{0}$, $u\geq 0$. For $u \in H^{1}_{0}$ which is not positve and general domain $\Omega$,  we use symmetrization arguments. Let $B_{R}$ be a ball having the same volume as $\Omega$ with $R=(|\Omega|/\omega_{n})^{1/n}$ and let $|u|^{*}$ be the symmetric decreasing rearrangement of the function $|u|$. Now note that for any $u \in H_{0}^{1}(\Omega)$, $|u|^{*} \in H_{0}^{1}(B_{R})$ and $|u|^{*}>0$. It is well known that the symmetrization does not change the $L^{p}$-norm, and that it decreases the Dirichlet energy, while increasing the integrals $\int_{\Omega}((\frac{n-2}{2})^{2}\frac{1}{|x|^{2}}+V(|x|)|u|^{2}dx$, since the weight $(\frac{n-2}{2})^{2}\frac{1}{|x|^{2}}+V(|x|)$  is a decreasing function of $|x|$.
Hence, (\ref{main-in}) holds for any $u \in H^{1}_{0}(\Omega)$. 

We shall need the following lemmas. 
 
\begin{lemma}\label{lemma}
Let $x(r)$ be a function in $C^{1}(0,R]$ that is a solution of
\begin{equation}\label{lem-ode}
rx'(r)+x^{2}(r)\leq -F(r), \ \ \ \ 0<r \leq R,
\end{equation}
where $F$ is a nonnegative continuous function. Then
\begin{equation}\label{zero}
\lim_{r\downarrow 0} x(r)=0.
\end{equation}
\end{lemma}
{\bf Proof:}  Divide equation (\ref{lem-ode}) by $r$ and integrate once. Then we have
\begin{equation}\label{lem-eq}
x(r)\geq \int_{r}^{R}\frac{|x(s)|^2}{s}ds+x(1)+\int_{r}^{R}\frac{F(s)}{s}ds.
\end{equation}
It follows  that $\lim_{r \downarrow 0}x(r)$ exists. In order to prove that this limit is zero, we claim that
\begin{eqnarray}\label{claim}
\int^{R}_{r}\frac{x^{2}(s)}{s}ds<\infty.
\end{eqnarray}
Indeed, otherwise we have
$G(r):=\int_{r}^{R}\frac{x^{2}(s)}{s}ds \rightarrow \infty$ as $r\rightarrow 0$.
From (\ref{lem-ode}) we have
\begin{equation*}
(-rG'(r))^{\frac{1}{2}}\geq G(r)+x(1)+\int^{R}_{r}\frac{F(s)}{s}ds.
\end{equation*}
Note that $F\geq 0$, and $G$ goes to infinity as $r$ goes to zero. Thus, for $r$ sufficiently small we have
$-rG'(r)\geq \frac{1}{2}G^{2}(r)$
hence,
$(\frac{1}{G(r)})'\geq\frac{1}{2}\ln (r)$,
which contradicts the fact that $G(r)$ goes to infinity as $r$ tends to zero. Thus, our claim is true and the limit  in (\ref{zero}) is indeed zero. \hfill $\Box$

\begin{lemma}\label{boundary.cond} If  the equation $\phi''+\frac{\phi'}{r}+v(r)\phi=0$ has a positive solution on some interval $(0, R)$, then  we have necessarily, 
\begin{equation}\label{boundary}
\hbox{$\liminf\limits_{r\rightarrow 0}r\frac{\varphi'(r)}{\varphi(r)}\geq 0$ \quad {\rm and} \quad  $\limsup\limits_{r\rightarrow R}\frac{\varphi'(r)}{\varphi(r)}<\infty$.}
\end{equation} 
\end{lemma}

{\bf Proof:} Since $\varphi(\delta)\geq 0$ and $\varphi(r)>0$ for $0<r<\delta$, it is obvious that $\varphi$ satisfies the second condition.  To obtain the first condition, set
$x(r)=r\frac{\varphi'(r)}{\varphi(r)}$. one may verify that 
 $x(r)$ satisfies the ODE:
\begin{equation*}
rx'(r)+x^{2}(r)= -F(r), \ \ \ \ {\rm for} \ \ 0<r\leq \delta,
\end{equation*}
where
$F(r)= r^2v(r)\geq 0.$
By Lemma \ref{lemma} we conclude that
$\lim_{r\downarrow 0}r\frac{\varphi'(r)}{\varphi(r)}=\lim_{r\downarrow 0}x(t)=0$. \hfill $\Box$ 

\begin{lemma} \label{super} Let $V$ be positive radial potential on the  ball $\Omega$ of radius $R$ in $R^{n}$ ($n \geq 3$). Assume that
\begin{eqnarray*}
\hbox{$\int_{\Omega}\left(|\nabla u|^{2}-(\frac{n-2}{2})^{2}\frac{|u|^{2}}{|x|^{2}}-V(|x|)|u|^{2}\right)dx\geq 0$ for all $u \in H_{0}^{1}(\Omega)$.}
\end{eqnarray*}
 Then there exists a $C^{2}$-supersolution to the equation
\begin{eqnarray}\label{pde}
-\Delta u-\left( \frac{n-2}{2}\right)^{2}\frac{u}{|x|^{2}}-V(|x|)u&=&0, \ \ \ \ {\rm in} \ \  \Omega, \\
u&>&0 \ \ \quad {\rm in} \ \  \Omega \setminus \{0\}, \\
 u&=&0 \ \quad {\rm in} \ \  \partial \Omega. 
\end{eqnarray}
 \end{lemma}
{\bf Proof:} Define
 \begin{eqnarray*}
\lambda_{1}(V):=\inf \{\frac{\int_{\Omega}|\nabla\psi|^{2}-(\frac{n-2}{2})^{2}|\psi|^{2}-V|\psi|^{2}}{\int_{\Omega}|\psi|^{2}}; \ \ \psi \in C^{\infty}_{0}(\Omega \setminus \{0\}) \}.
\end{eqnarray*}
By our assumption $\lambda(V)\geq 0$. Let $(\phi_{n}, \lambda^{n}_{1})$ be the first eigenpair for the problem
\begin{eqnarray*}
(L-\lambda_{1}(V)-\lambda^{r}_{1})\phi_{r}&=&0 \ \ on \ \ \Omega \setminus B_{\frac{R}{n}}\\
\phi(r)&=&0 \ \ on\ \ \partial (\Omega \setminus B_{\frac{R}{n}}),
\end{eqnarray*}
where $L=-\Delta-(\frac{n-2}{2})^{2}\frac{1}{|x|^{2}}-V$, and $B_{\frac{R}{n}}$ is a ball of radius $\frac{R}{n}$, $n\geq 2$ . The eigenfunctions can be chosen in such a way that $\phi_{n}>0$ on $\Omega \setminus B_{\frac{R}{n}}$ and $\varphi_{n}(b)=1$, for some $b \in \Omega$ with $\frac{R}{2}<|b|<R$.

Note that $\lambda^{n}_{1}\downarrow 0$ as $n \rightarrow \infty$. Harnak's inequality yields that for any compact subset $K$, $\frac{max_{K}\phi_{n}}{max_{K}\phi_{n}}\leq C(K)$ with the later constant being independant of $\phi_{n}$.  Also standard elliptic estimates also yields that the family $(\phi_{n})$ have also uniformly bounded derivatives on compact sets $\Omega-B_{\frac{R}{n}}$.  \\
Therefore, there exists a subsequence $(\varphi_{n_{l_{2}}})_{l_{2}}$ of ($\varphi_{n})_{n}$ such that $(\varphi_{n_{l_{2}}})_{l_{2}}$ converges to some $\varphi_{2} \in C^{2}(\Omega \setminus B(\frac{R}{2}))$.  Now consider  $(\varphi_{n_{l_{2}}})_{l_{2}}$ on $\Omega \setminus B(\frac{R}{3})$.  Again there exists a subsequence  $(\varphi_{n_{l_{3}}})_{l_{3}}$ of $(\varphi_{n_{l_{2}}})_{l_{2}}$ which converges to $\varphi_{3} \in C^{2}( \Omega \setminus B(\frac{R}{3}))$,  and $\varphi_{3}(x)=\varphi_{2}(x)$ for all $x \in \Omega \setminus B(\frac{R}{2})$. By repeating this argument we get a supersolution $\varphi \in C^{2}( \Omega \setminus\{ 0\})$ i.e. $L\varphi \geq 0$, such that $\varphi>0$ on $\Omega \setminus \{0\}$. \hfill $\square$

\begin{lemma} \label{OSI-LEM} Let $a$ be a locally integrable function on $\R$, then the following statements are equivalent. 
\begin{enumerate}
\item  $z''(s)+a(s)z(s)=0$, has a strictly positive solution on $(b,\infty)$.

\item There exists a function $\psi \in C^{1}(b,\infty)$ such that 
$\psi'(r)+\psi^{2}(r)+a(t)\leq 0, \ \ for \ \ r>b.$
\end{enumerate}
Consequently,  the equation  $y''+\frac{1}{r}y'+v(r)y= 0$ has a positive supersolution on $(0, \delta)$ if and 
only  if it has a positive solution on $(0, \delta)$. 
 
 \end{lemma}
 {\bf Proof:} That 1) and 2) are equivalent follows from the work of Wintner \cite{win1, win2}, a proof of which may be found in \cite{Har}). \\
To prove the rest, we note that  the change of variable 
 $z(s)=y(e^{-s})$ maps the equation $y''+\frac{1}{r}y'+v(r)y=0$ into $z''+e^{-2s}v(e^{-s})z(s)=0$. On the other hand, the change of variables $\psi (t)=\frac{-e^{-t}y'(e^{-t})}{y(e^{-t})}$ maps  $y''+\frac{1}{r}y'+v(r)y$ into $\psi'(t)+\psi^{2}(t)+e^{-2t}v(e^{-t})$. This proves the lemma. \hfill $\square$\\
 
  {\bf Proof of Theorem \ref{main}:} The implication 1) implies 2) follows immediately from Proposition \ref{main-pro} and Lemma \ref{boundary.cond}. It is valid for any smooth bounded domain provided $v$ is assumed to be non-decreasing on $(0, R)$. this condition is not needed if the domain is a ball of radius $R$.
  
  To show that 2) implies 1), we assume that inequality (${\rm H}_{V}$) holds on a ball $\Omega$ of radius $R$, and then apply Lemma (\ref{super}) to obtain a $C^{2}$-supersolution for the equation (\ref{pde}). Now take the surface average of $u$, that is
\begin{equation}
w(r)=\frac{1}{n\omega_{w} r^{n-1}}\int_{\partial B_{r}} u(x)dS=\frac{1}{n\omega_{n}} \int_{|\omega|=1}u(r\omega)d\omega >0,
\end{equation}
where $\omega_{n}$ denotes the volume of the unit ball in $R^{n}$. We may assume that the unit ball is contained in $\Omega$ (otherwise we just use a smaller ball). By a standard calculation we get
\begin{equation}
w''(r)+\frac{n-1}{r}w'(r)\leq \frac{1}{n\omega_{n}r^{n-1}}\int_{\partial B_{r}}\Delta u(x)dS.
\end{equation}
Since $u(x)$  is a supersolution of (\ref{pde}), $w$ satisfies the inequality:
\begin{equation}\label{ode}
w''(r)+\frac{n-1}{r}w'(r)+(\frac{n-2}{2})^{2}\frac{w(r)}{r^{2}}\leq -v(r)w(r), \ \ \ \ for\ \ \ 0<r<R.
\end{equation}
Now define
\begin{equation}
\varphi(r)=r^{\frac{n-2}{2}}w(r), \ \ \ \ in\ \ 0<r<R.
\end{equation}
Using (\ref{ode}), a straightforward calculation shows that $\varphi$ satisfies the following inequality
\begin{equation}
\varphi''(r)+\frac{\varphi'(r)}{r}\leq -\varphi(r)v(r), \ \ \ \ for \ \ 0<r<R.
\end{equation}
By  Lemma \ref{OSI-LEM} we may conclude that the equation $y''(r)+\frac{1}{r}y' +v(r)y=0$ has actually a positive solution $\phi$ on $(0, R)$.

To establish formula (\ref{best.constant}), it is clear that by the sufficient condition $c(V) \geq c$ whenever 
$y''(r)+\frac{1}{r}y' +cv(r)y=0$ has  a positive solution on $(0, R)$. On the other hand, the necessary condition yields that $y''(r)+\frac{1}{r}y' +c(V)v(r)y=0$ has a positive solution on $(0, R)$. The proof is now complete. \hfill $\Box$


\section{Applications}

In this section we start by applying Theorem \ref{main} to recover in a relatively simple and unified way,  all previously known improvements of Hardy's inequality.  For that  we need to investigate whether the ordinary differential equation
\begin{equation}
y''+\frac{y'}{r}+v(r)y(r)=0,
\end{equation}
corresponding to a potential $v$ has a positive solution $\phi$ on $(0,\delta)$ for some $\delta>0$.  In this case,  $\psi (r)=\phi (\frac{\delta r}{R})$ is a solution for $y''(r)+\frac{1}{r}y' +\frac{\delta^2}{R^2}v(\frac{\delta}{R}r)y=0$ on $(0, R)$, which means that the scaled potential $V_\delta (x)=\frac{\delta^2}{R^2}V(\frac{\delta}{R}x)$ yields an improved Hardy formula (${\rm H}_{V_\delta}$) on a ball of radius $R$, with constant larger than one. 
 Here is an immediate application of this criterium.\\

{\bf 1) The Brezis-V\'{a}zquez improvement  \cite{BV}:} Here we need to show that  we can have an improved inequality with a constant potential. In this case,  the best constant for which the equation
\begin{equation}\label{BV}
y''+\frac{y'}{r}+cy(r)=0,
\end{equation}
 has a positive solution on $(0,R)$, with $R=(| \Omega | / \omega_{n})^{\frac{1}{n}}$ is $z^{2}_{0}\omega^{2/n}_{n}|\Omega|^{-2/n}$. Indeed, if $z_0$ is the first root of the solution of the Bessel equation $y''+\frac{y'}{r}+y(r)=0$, then the solution of (\ref{BV}) in this case is  the Bessel function $\varphi(r)=J_{0}(\frac{rz_{0}}{R})$.  This readily gives the result of  Brezis-V\'{a}zquez  mentioned in the introduction. \\

{\bf 2) The Adimurthi et al. improvement  \cite{ACR}:} In this case, one easily sees that the functions 
$\varphi_{j}(r)=( \prod^{j}_{i=1}log^{(i)}\frac{\rho}{r})^{\frac{1}{2}}$ is a solution of the equation 
\begin{equation*}
-\frac{\varphi_{j}'(r)+r\varphi_{j}''(r)}{r\varphi_{j}(r)}=\frac{1}{4r^2}(\prod^{n}_{i=1}\log^{(i)}\frac{\rho}{r})^{-2},
\end{equation*}
on $(0, R)$, which means that the inequality $({\rm H}_V)$ holds for the potential $V(x)=\frac{1}{4|x|^2}(\prod^{n}_{i=1}\log^{(i)}\frac{\rho}{|x|})^{-2}$ which yields the result of Adimurthi et al. In the following, we use our characterization to show that the constant appearing in the above improvement is indeed the best constant in the following sense:\\
\begin{equation}
\frac{1}{4}=\inf_{u \in H_{0}^{1}(\Omega)\setminus \{0\}}\frac{\int_{\Omega}|\nabla u|^{2}dx - ( \frac{n-2}{2})^{2} \int_{\Omega}\frac{|u|^{2}}{|x|^{2}}dx-\frac{1}{4}\sum^{m-1}_{j=1}\int_{\Omega}\frac{|u|^{2}}{|x|^{2}}\big( \prod^{j}_{i=1}log^{(i)}\frac{\rho}{|x|}\big)^{-2}}{\int_{\Omega}\frac{|u|^{2}}{|x|^{2}}\big( \prod^{m}_{i=1}log^{(i)}\frac{R}{|x|}\big)^{-2}},
\end{equation}
fo all $1 \leq m \leq k $. We proceed by contradiction, and assume that 
\begin{equation*}
\frac{1}{4}+\lambda=\inf_{u \in H_{0}^{1}(\Omega) \setminus \{0\}}\frac{\int_{\Omega}|\nabla u |^{2}dx - ( \frac{n-2}{2})^{2} \int_{\Omega}\frac{|u|^{2}}{|x|^{2}}dx-\frac{1}{4}\sum^{m-1}_{j=1}\int_{\Omega}\frac{|u|^{2}}{|x|^{2}}\big( \prod^{j}_{i=1}log^{(i)}\frac{\rho}{|x|}\big)^{-2}}{\int_{\Omega}\frac{|u|^{2}}{|x|^{2}}\big( \prod^{m}_{i=1}log^{(i)}\frac{\rho}{|x|}\big)^{-2}},
\end{equation*}
and $\lambda>0$. From Theorem \ref{main} we deduce that there exists a positive function $\varphi$ such that
\begin{equation*}
-\frac{\varphi'(r)+r\varphi''(r)}{\varphi(r)}=\frac{1}{4}\sum^{m-1}_{j=1}\frac{1}{r}\big( \prod^{j}_{i=1}log^{(i)}\frac{\rho}{r}\big)^{-2}+(\frac{1}{4}+\lambda)\frac{1}{r}\big( \prod^{m}_{i=1}log^{(i)}\frac{\rho}{r}\big)^{-2}.
\end{equation*}
 Now define $f(r)=\frac{\varphi(r)}{\varphi_{m}(r)}>0$, and calculate, 
\[\frac{\varphi'(r)+r\varphi''(r)}{\varphi(r)}=\frac{\varphi_{m}'(r)+r\varphi_{m}''(r)}{\varphi_{m}(r)}+\frac{f'(r)+rf''(r)}{f(r)}-\frac{f'(r)}{f(r)}\sum_{i=1}^{m}\frac{1}{\prod^{i}_{j=1}\log^{j}(\frac{\rho}{r})}.\]
Thus,
\begin{equation}\label{main-rel}
\frac{f'(r)+rf''(r)}{f(r)}-\frac{f'(r)}{f(r)}\sum_{i=1}^{m}\frac{1}{\prod^{i}_{j=1}\log^{j}(\frac{\rho}{r})}=-\lambda\frac{1}{r}\big( \prod^{m}_{i=1}log^{(i)}\frac{\rho}{r}\big)^{-2}.
\end{equation}
If now $f'(\alpha_{n})=0$ for some sequence $\{\alpha_{n}\}^{\infty}_{n=1}$ that converges to zero, then there exists a sequence $\{\beta_{n}\}^{\infty}_{n=1}$ that also converges to zero, such that  $f''(\beta_{n})=0$, and $f'(\beta_{n})>0$. But this contradicts (\ref{main-rel}), which means that $f$ is eventually monotone for $r$ small enough. We consider the two cases according to whether $f$ is increasing or decreasing:\\ \\
Case I:  Assume $f'(r)>0$ for $r>0$ sufficiently small. Then we will have
\[\frac{(rf'(r))'}{rf'(r)}\leq \sum_{i=1}^{m}\frac{1}{r\prod^{i}_{j=1}\log^{j}(\frac{\rho}{r})}.\]
Integrating once we get
\[f'(r)\geq \frac{c}{r\prod^{m}_{j=1}\log^{j}(\frac{\rho}{r})},\]
for some $c>0$.
Hence, $\lim_{r\rightarrow 0}f(r)=-\infty$ which is a contradiction. \\ \\
Case II: Assume $f'(r)<0$ for $r>0$ sufficiently small. Then
\[\frac{(rf'(r))'}{rf'(r)}\geq \sum_{i=1}^{m}\frac{1}{r\prod^{i}_{j=1}\log^{j}(\frac{\rho}{r})}. \]
Thus,
\begin{equation}\label{estim}
f'(r)\geq- \frac{c}{r\prod^{m}_{j=1}\log^{j}(\frac{\rho}{r})},
\end{equation}
for some $c>0$ and $r>0$ sufficiently small. On the other hand
\[\frac{f'(r)+rf''(r)}{f(r)}\leq-\lambda \sum^{m}_{j=1} \frac{1}{r}\big( \prod^{j}_{i=1}log^{(i)}\frac{R}{r}\big)^{-2}\leq -\lambda(\frac{1}{\prod^{m}_{j=1}\log^{j}(\frac{\rho}{r})})'.\]
Since $f'(r)<0$, there exists $l$ such that $f(r)>l>0$ for $r>0$ sufficiently small. From the above inequality we then have
\[bf'(b)-af'(a)<-\lambda l(\frac{1}{\prod^{m}_{j=1}\log^{j}(\frac{\rho}{b})}-\frac{1}{\prod^{m}_{j=1}\log^{j}(\frac{\rho}{a})}).\]
From (\ref{estim}) we have $\lim_{a \rightarrow 0}af'(a)=0$. Hence, 
\[bf'(b)<-\frac{\lambda l }{\prod^{m}_{j=1}\log^{j}(\frac{\rho}{b})},\]
for every $b>0$, and
\[f'(r)<- \frac{\lambda l }{r\prod^{m}_{j=1}\log^{j}(\frac{\rho}{r})},\]
for $r>0$ sufficiently small. Therefore, 
\[\lim_{r\rightarrow 0}f(r)=+\infty,\]
and by choosing $l$ large enouph (e.g.,  $l>\frac{c}{\lambda})$ we get to contradict  $(\ref{estim})$ and the proof is now complete. \\ \\

{\bf 3) The Filippas and Tertikas improvement \cite{FT}:} Let $D\geq \sup_{x \in \Omega}|x|$, and define
\begin{eqnarray*}
\varphi_{k}(r)=(X_{1}(\frac{r}{D})X_{2}(\frac{r}{D}) \ldots X_{i-1}(\frac{r}{D})X_{i}(\frac{r}{D}))^{-\frac{1}{2}}, \ \ \ \ i=1,2, \ldots .
\end{eqnarray*}
Using the fact that
$X'_{k}(r)=\frac{1}{r}X_{1}(r)X_{2}(r) \ldots X_{k-1}(r)X^{2}_{k}(r)$ for  $k=1,2, \ldots $, 
we get
\begin{eqnarray*}
-\frac{\varphi_{k}'(r)+r\varphi_{k}''(r)}{\varphi_{k}(r)}=\frac{1}{4r}X^{2}_{1}(\frac{r}{D})X^{2}_{2}(\frac{r}{D}) \ldots X^{2}_{k-1}(\frac{r}{D})X^{2}_{k}(\frac{r}{D}).
\end{eqnarray*}
This means that the inequality $({\rm H}_V)$ holds for the potential $V(x)=\frac{1}{4|x|^{2}}X^{2}_{1}(\frac{|x|}{D})X^{2}_{2}(\frac{|x|}{D}) \ldots X^{2}_{k-1}(\frac{|x|}{D})X^{2}_{k}(\frac{|x|}{D})$, which yields the result  of Filippas and Tertikas \cite{FT}. We now identify the best constant by showing that:
  \begin{equation}
\frac{1}{4}=\inf_{u \in H_{0}^{1}(\Omega) \setminus \{0\}}\frac{\int_{\Omega}|\nabla u |^{2}dx - ( \frac{n-2}{2})^{2} \int_{\Omega}\frac{|u|^{2}}{|x|^{2}}dx-\frac{1}{4}\sum^{m-1}_{j=1}\int_{\Omega}\frac{|u|^{2}}{|x|^{2}}X^{2}_{1}(\frac{|x|}{D})X^{2}_{2}(\frac{|x|}{D}) \ldots X^{2}_{j-1}(\frac{|x|}{D})X^{2}_{j}(\frac{|x|}{D})}{\int_{\Omega}\frac{|u|^{2}}{|x|^{2}}X^{2}_{1}(\frac{|x|}{D})X^{2}_{2}(\frac{|x|}{D}) \ldots X^{2}_{m-1}(\frac{|x|}{D})X^{2}_{m}(\frac{|x|}{D})},
\end{equation}
fo all $1 \leq m \leq k $. We proceed again by contradiction and in a way very similar to the above case. Indeed, assuming  that 
\begin{equation*}
\frac{1}{4}+\lambda=\inf_{u \in H_{0}^{1}(\Omega) \setminus \{0\}}\frac{\int_{\Omega}|\nabla u |^{2}dx - ( \frac{n-2}{2})^{2} \int_{\Omega}\frac{|u|^{2}}{|x|^{2}}dx-\frac{1}{4}\sum^{m-1}_{j=1}\int_{\Omega}\frac{|u|^{2}}{|x|^{2}}X^{2}_{1}(\frac{|x|}{D})X^{2}_{2}(\frac{|x|}{D}) \ldots X^{2}_{j-1}(\frac{|x|}{D})X^{2}_{j}(\frac{|x|}{D})}{\int_{\Omega}\frac{|u|^{2}}{|x|^{2}}X^{2}_{1}(\frac{|x|}{D})X^{2}_{2}(\frac{|x|}{D}) \ldots X^{2}_{m-1}(\frac{|x|}{D})X^{2}_{m}(\frac{|x|}{D})},
\end{equation*}
and $\lambda>0$, we use again Theorem \ref{main} to find a positive function $\varphi$ such that
\begin{eqnarray*}
-\frac{\varphi'(r)+r\varphi''(r)}{\varphi(r)}&=&\frac{1}{4}\sum^{m-1}_{j=1}\frac{1}{r}X^{2}_{1}(\frac{r}{D})X^{2}_{2}(\frac{r}{D}) \ldots X^{2}_{j-1}(\frac{r}{D})X^{2}_{j}(\frac{r}{D})\\
&+&(\frac{1}{4}+\lambda)\frac{1}{r}X^{2}_{1}(\frac{r}{D})X^{2}_{2}(\frac{r}{D}) \ldots X^{2}_{m-1}(\frac{r}{D})X^{2}_{m}(\frac{r}{D}).
\end{eqnarray*}
Setting $f(r)=\frac{\varphi(r)}{\varphi_{m}(r)}>0$, we have
\[\frac{\varphi'(r)+r\varphi''(r)}{\varphi(r)}=\frac{\varphi_{m}'(r)+r\varphi_{m}''(r)}{\varphi_{m}(r)}+\frac{f'(r)+rf''(r)}{f(r)}-\frac{f'(r)}{f(r)}\sum_{i=1}^{m}\prod_{j=1}^{i}X_{j}(\frac{r}{D}).\]
Thus,
\begin{equation}\label{main-rel2}
\frac{f'(r)+rf''(r)}{f(r)}-\frac{f'(r)}{f(r)}\sum_{i=1}^{m}\prod_{j=1}^{i}X_{j}(\frac{r}{D})=-\lambda\frac{1}{r}\prod_{j=1}^{m}X^{2}_{j}(\frac{r}{D}).
\end{equation}
Arguing as before, we deduce that $f$ is eventually monotone for $r$ small enough,  and we consider two cases:\\ \\
Case I:  If $f'(r)>0$ for $r>0$ sufficiently small, then we will have
\[\frac{(rf'(r))'}{rf'(r)}\leq \sum_{i=1}^{m}\frac{1}{r}\prod_{j=1}^{i}X_{j}(\frac{r}{D}).\]
Integrating once we get
\[f'(r)\geq \frac{c}{r}\prod_{j=1}^{m}X_{j}(\frac{r}{D}),\]
for some $c>0$, and 
therefore $\lim_{r\rightarrow 0}f(r)=-\infty$ which is a contradiction. \\ \\
Case II: Assume $f'(r)<0$ for $r>0$ sufficiently small. Then
\[\frac{(rf'(r))'}{rf'(r)}\geq \sum_{i=1}^{m}\frac{1}{r}\prod_{j=1}^{i}X_{j}(\frac{r}{D})\]
Thus,
\begin{equation}\label{estim2}
f'(r)\geq- \frac{c}{r}\prod_{j=1}^{m}X_{j}(\frac{r}{D}),
\end{equation}
for some $c>0$ and $r>0$ sufficiently small. On the other hand
\[\frac{f'(r)+rf''(r)}{f(r)}\leq-\lambda \sum^{m}_{j=1} \frac{1}{r}\prod^{j}_{i=1}X^{2}_{j}\leq -\lambda(\prod^{m}_{j=1}X_{j}(\frac{r}{D})  )'.\]
Since $f'(r)<0$, we may  assume $f(r)>l>0$ for $r>0$ sufficiently small, and from the above  inequality we have
\[bf'(b)-af'(a)<-\lambda l( \prod^{m}_{j=1}X_{j}(\frac{b}{D})-  \prod^{m}_{j=1}X_{j}(\frac{a}{D})).\]
From (\ref{estim2}) we have $\lim_{a \rightarrow 0}a f'(a)=0$. Hence, 
\[f'(r))<-\frac{\lambda l}{r}\prod^{m}_{j=1}X_{j}(\frac{r}{D}),\]
for $r>0$ sufficiently small. Therefore, 
\[\lim_{r\rightarrow 0}f(r)=+\infty,\]
and by choosing  $l$ large enouph (i.e. $l>\frac{c}{\lambda})$ we contradict $(\ref{estim2})$ and the proof is complete. \hfill $\square$\\

We shall now make the connection between improved Hardy inequalities and 
the existence of non-oscillatory solutions (i.e., those $z(s)$ such that $z(s)>0$ for $s>0$ sufficiently large) for the second order linear differential equations
\begin{equation}\label{OSI_ODE}
z''(s)+a(s)z(s)=0.
\end{equation}
Interesting results in this direction were established by many authors (see \cite{Har, Hua, win1,win2, won}). Here is a typical  criterium about the oscillatory properties of equation (\ref{OSI_ODE}):
 \begin{enumerate}
\item   If $\hbox{$\limsup_{t\rightarrow \infty } t\int^{\infty}_{t}a(s)ds<\frac{1}{4}$, } $ then Eq. (\ref{OSI_ODE}) is non-oscillatory. 
\item  If 
$\hbox{$\liminf_{t\rightarrow \infty } t\int^{\infty}_{t}a(s)ds>\frac{1}{4}$,}$
 then Eq. (\ref{OSI_ODE}) is oscillatory.
\end{enumerate}
This result combined with Theorem \ref{main} and Lemma \ref{OSI-LEM} clearly yields Corollary \ref{osc}. \\
 
{\bf Proof of Corollary \ref{dual}:} It follows from H\"older's inequality that
\[(\int_{\Omega}V(|x|)u^{2}(x)dx)^{\frac{1}{s}}\geq \frac{\int_{\Omega}u^{\frac{2}{s}}(x)dx}{(\int_{\Omega}V^{-\frac{r}{s}}(|x|)dx)^{\frac{1}{r}}},\]
where $s\geq1$ and $\frac{ 1}{s}+\frac{1}{r}=1$.  Letting $p=\frac{2}{s}$, we get 
\[\int_{\Omega}V(|x|)u^{2}(x)dx\geq (\int_{\Omega}u^{p}(x)dx)^{\frac{2}{p}}\frac{1}{||V^{-1}(|x|)||_{L^{\frac{p}{2-p}}(\Omega)}}.\]
Inequality (\ref{dual}) now follows from Theorem \ref{main}. \hfill $\square$\\ \\
{\bf Proof of Corollary \ref{eigenvalue}:} Define the functional 
\begin{equation}
F_{\mu}(u)=\int_{\Omega}|\nabla u(x)|^{2}dx-\mu \int_{\Omega}\frac{u^{2}(x)}{|x|^{2}}dx,
\end{equation}
which is continuous, Gateaux differentiable and coercive on $H^{1}_{0}(\Omega)$. Let  $u_{\mu}>0$ be a  minimizer of $F_\mu$ over the manifold $M=\{u \in H_{0}^{1}(\Omega) | \ \ \int_{\Omega} u^{2}(x)V(x)=1\}$ and assume $\lambda^{1}_{\mu}$ is the infimum. It is clear that $\lambda^{1}_{\mu}>0$. By standard arguments we can conclude that  $u_{\mu}$ is a weak solution of $(E_{V,\mu})$. The rest of the proof follows from Corollary \ref{osc} and the fact that
\[\lambda_{1}(V)=\lim_{\mu \rightarrow \mu_{n}}\lambda^{1}_{\mu}=\inf_{u \in H_{0}^{1}(\Omega)\setminus \{0 \}}\frac{\int_{\Omega}(|\nabla u|^{2}-\mu_{n}\frac{u^{2}(x)}{|x|^{2}})dx}{\int_{\Omega}|u(x)|^{2}V(x)dx}.\]

 \end{document}